\documentclass[12pt,twoside]{article}

 \usepackage{float}
\usepackage{graphicx}
\usepackage{epstopdf}
\usepackage{graphicx}
\usepackage{epic}
\usepackage{multirow}
\usepackage{tikz}
\usepackage{xcolor}
\usepackage{threeparttable}
\usetikzlibrary{arrows,shapes,chains}

\usepackage{graphicx}
\usepackage{makecell}
\renewcommand{\paragraph}{\roman{paragraph}}
\usepackage[a4paper]{geometry}
\makeatletter
\renewcommand\title[1]{\gdef\@title{\reset@font\Large\bfseries #1}}
\renewcommand\section{\@startsection {section}{1}{\z@}%
                                   {-3.5ex \@plus -1ex \@minus -.2ex}%
                                   {2.3ex \@plus.2ex}%
                                   {\normalfont\large\bfseries}}
\renewcommand\subsection{\@startsection{subsection}{2}{\z@}%
                                     {-3ex\@plus -1ex \@minus -.2ex}%
                                     {1.5ex \@plus .2ex}%
                                     {\normalfont\normalsize\bfseries}}
\renewcommand\subsubsection{\@startsection{subsubsection}{3}{\z@}%
                                     {-2.5ex\@plus -1ex \@minus -.2ex}%
                                     {1.5ex \@plus .2ex}%
                                     {\normalfont\normalsize\bfseries}}

\def\@runningauthor{}\newcommand{\runningauthor}[1]{\def\runningauthor{#1}}
\def\@runningtitle{}\newcommand{\runningtitle}[1]{\def\runningtitle{#1}}

\renewcommand{\ps@plain}{%
\renewcommand{\@evenhead}{\footnotesize\scshape \hfill\runningauthor\hfill}
\renewcommand{\@oddhead}{\footnotesize\scshape \hfill\runningtitle\hfill}}

\newcommand{\bm}[1]{\mbox{\boldmath{$#1$}}}
\g@addto@macro\bfseries{\boldmath}

\makeatother

\usepackage{amsthm,amsmath,amssymb}
\usepackage{cite}
\usepackage{graphicx}

\usepackage[colorlinks=true,citecolor=black,linkcolor=black,urlcolor=blue]{hyperref}

\theoremstyle{plain}

\theoremstyle{definition}

\theoremstyle{remark}

\runningauthor{}

\date{}

\begin{document}
\title{The von Bahr-Esseen type inequality under sub-linear
expectations and applications\thanks{Supported
by the National Natural Science Foundation of China (12201079), the Natural Science Foundation of Anhui
Province (2108085MA06, 2108085QA15), and the Excellent Scientific Research and Innovation Team of Anhui Colleges (2022AH010098).}}
\author{Yi Wu$^{a}$,~~~~Xuejun Wang$^{b}$\thanks{Corresponding author. Email address: wxjahdx2000@126.com}\\
{\small   $^{a}$ School of Big Data and Artificial Intelligence, Chizhou University, Chizhou, 247000, P.R. China}\\
{\small   $^{b}$ School of Big Data and Statistics, Anhui University, Hefei, 230601, P.R. China}}
\date{}
    \maketitle
\begin{abstract}
Moment inequalities play important roles in probability limit theory and mathematical statistics. In this work, the von Bahr-Esseen type inequality for extended negatively dependent random variables under sub-linear expectations is established successfully. By virtue of the inequality, we further obtain the Kolmogorov type weak law of large numbers for partial sums and the complete convergence for weighted sums, which extend and improve corresponding results in sub-linear expectation space.
\end{abstract}
\vspace{2mm}


\section{Introduction}

It is known that the moment inequalities are very powerful tools in classical probability theory and statistics. Due to these inequalities, many elegant results in probability theory and statistical applications were established successively in the past decades. For a sequence $\{X_{n},n\geq1\}$ of independent and identically distributed random variables with mean zero and $E|X_{n}|^{p}<\infty$ for some $1\leq p\leq2$, von Bahr and Esseen (1965) proved that, for each $n\geq1$,
\begin{equation}
E\left|\sum_{i=1}^{n}X_{i}\right|^{p}\leq C_{p}\sum_{i=1}^{n}E|X_{i}|^{p},\label{I1}\end{equation}
where $C_{p}$ is a positive constant depending only on $p$. The aforementioned inequality is called the $p$-th von Bahr-Esseen moment inequality. In view of the significance of this powerful tool, many scholars generalized it to various dependence settings. For example, Chatterji (1969) showed that \eqref{I1} holds for martingale difference sequence; Bryc and Smolenski (1993) proved it for weakly correlated random variables; Shao (2000) extended it to negatively associated random variables; Asadian et al. (2006) established \eqref{I1} for negatively orthant dependent random variables; Wang et al. (2014) generalized it to widely orthant dependent random variables; Chen et al. (2014) innovatively proved \eqref{I1} with $1<p<2$ by that with $p=2$ for pairwise independent as well as pairwise negatively quadrant dependent random variables. Later on, Wu et al. (2021) obtained \eqref{I1} for asymptotic negatively associated random variables as well as $m$-asymptotic negatively associated random variables by adopting the method in Chen et al. (2014).

Most of moment inequalities including \eqref{I1} are based on the linearity of expectations and probability measures. However, such an additivity assumption is not feasible in many areas of applications because many uncertain phenomena cannot be well modelled by using additive probabilities or additive expectations.

In recent years, the sub-linear expectations receive a lot of attention from statisticians. There are many entitative differences between the classical probability space and the sub-linear expectation space. For example, in the classical probability space, the probability measures and mathematical expectations are always assumed to be additive. The linear expectation and determinant statistics are built on the distribution-certainty
or model-certainty. However, in practical fields such as risk measure and super-hedging in finance, this assumption is not always plausible. For related
references, we refer the readers to El Karoui et al. (1997), Artzner et al. (1997), Chen and Epstein (2002), and F\"{o}llmer and Schied
(2004) among others. Without distribution-certainty, the resulting expectation is usually nonlinear. Therefore, the theory and methodology
of nonlinear expectation have been well developed and received much attention
in some application fields such as finance risk measure and control in the past decades. A typical example of the nonlinear expectation, called $g$-expectation,
was introduced by Peng (1997) in the framework of backward stochastic differential
equations. Peng (2006) proposed the $G$-expectation and its
related versions as a further development. For more details about the applications of the nonlinear expectation, one can refer to Denis and Martini (2006), Gilboa (1987), Marinacci (1999), and Peng (1999, 2008a) for instance. Peng (2006, 2008a, 2008b) firstly introduced a reasonable framework of the sub-linear expectation in a general function space by relaxing the linear property of the classical expectation to the sub-additivity and positive homogeneity.

Let $(\Omega,\mathcal{F})$ be a given measurable space and $\mathcal{H}$ be a linear space of real functions defined on $(\Omega,\mathcal{F})$ such that if $X_{1},\cdots,X_{n}\in \mathcal{H}$ then $\varphi(X_{1},\cdots,X_{n})\in \mathcal{H}$ for each $\varphi\in C_{l,Lip}(\mathbb{R}^{n})$, where $C_{l,Lip}(\mathbb{R}^{n})$ denotes the linear space of (local Lipschitz) functions $\varphi$ satisfying
\begin{eqnarray*}|\varphi(\bm{x})-\varphi(\bm{y})|\leq C(1+|\bm{x}|^{m}+|\bm{y}|^{m})|\bm{x}-\bm{y}|,\textrm{ for any }\bm{x},\bm{y}\in\mathbb{R}^{n},\end{eqnarray*}
for some $C>0$, $m\in\mathbb{N}$ depending on $\varphi$. $\mathcal{H}$ can be considered as a space of ``random variables". In this case we denote $X\in\mathcal{H}$. We also denote $C_{b,Lip}(\mathbb{R}^{n})$ to be the space of  bounded Lipschitz functions $\varphi$ satisfying
\begin{eqnarray*}|\varphi(\bm{x})|\leq C,~~|\varphi(\bm{x})-\varphi(\bm{y})|\leq C|\bm{x}-\bm{y}|,\textrm{ for any }\bm{x},\bm{y}\in\mathbb{R}^{n},\end{eqnarray*}
for some $C>0$ depending on $\varphi$.

\noindent\textbf{Definition 1.1.}~~{\it A sub-linear expectation $\hat{\mathbb{E}}$ on $\mathcal{H}$ is a function $\hat{\mathbb{E}}: \mathcal{H}\rightarrow\bar{\mathbb{R}}$ satisfying the following properties: for all $X,Y\in\mathcal{H}$, we have

\noindent(a) Monotonicity: If $X\geq Y$, then $\hat{\mathbb{E}}[X]\geq\hat{\mathbb{E}}[Y]$;

\noindent(b) Constant preserving: $\hat{\mathbb{E}}[c]=c$, where $c$ is a constant;

\noindent(c) Sub-additivity: $\hat{\mathbb{E}}[X+Y]\leq\hat{\mathbb{E}}[X]+\hat{\mathbb{E}}[Y]$ whenever $\hat{\mathbb{E}}[X]+\hat{\mathbb{E}}[Y]$ is not the form $+\infty-\infty$ or $-\infty+\infty$;

\noindent(d) Positive homogeneity: $\hat{\mathbb{E}}[\lambda X]=\lambda\hat{\mathbb{E}}[X]$, $\lambda>0$.}

 Here $\bar{\mathbb{R}}=[-\infty,+\infty]$. The triple $(\Omega,\mathcal{H},\hat{\mathbb{E}})$ is called a sub-linear expectation space. For a given
sub-linear expectation $\hat{\mathbb{E}}$, the conjugate expectation \text{\Large $\hat{\varepsilon}$} of $\hat{\mathbb{E}}$ is denoted by
\begin{eqnarray*}\text{\Large $\hat{\varepsilon}$}[X]:=-\hat{\mathbb{E}}[-X],\textrm{ for any }X\in \mathcal{H}.\end{eqnarray*}
It is easy to verify that \text{\Large $\hat{\varepsilon}$}$[X]\leq\hat{\mathbb{E}}[X]$, $\hat{\mathbb{E}}[X+c]=\hat{\mathbb{E}}[X]+c$ and $\hat{\mathbb{E}}[X-Y]\geq\hat{\mathbb{E}}[X]-\hat{\mathbb{E}}[Y]$ for all $X,Y\in\mathcal{H}$ with $\hat{\mathbb{E}}[Y]$ being finite. Furthermore, if $\hat{\mathbb{E}}[|X|]$ is finite, then both \text{\Large $\hat{\varepsilon}$}$[X]$ and $\hat{\mathbb{E}}[X]$ are finite.

Let $(\Omega,\mathcal{H},\hat{\mathbb{E}})$ be a sub-linear expectation space. Peng (2010) initiated and proved that there exists a family of finitely additive probabilities $\{P_{\theta}\}_{\theta\in\Theta}$ defined on $(\mathbb{R},\mathcal{B}(\mathbb{R}))$ such that for each $X\in\mathcal{H}$,
\begin{eqnarray*}\hat{\mathbb{E}}[X]=\sup_{\theta\in\Theta}E_{P_{\theta}}[X].\end{eqnarray*}
Let $\mathcal{P}$ be a family of probability measures defined on the measurable space ($\Omega,\mathcal{F}$). Define an upper expectation $\hat{\mathbb{E}}[\cdot]$ by
$\hat{\mathbb{E}}[X]=\sup_{P\in\mathcal{P}}E_{P}[X],$
where $E_{P}$ denotes the classical expectation under probability $P$.
Then $\hat{\mathbb{E}}$ is a sub-linear expectation.

Peng (2008b) put forward the definitions of independence and identical distribution through the sub-linear expectation. {\it A random vector $\bm{Y}=(Y_{1},\cdots,Y_{n})$, $Y_{i}\in \mathcal{H}$ is said to be independent to another random vector $\bm{X}=(X_{1},\cdots,X_{n})$, $X_{i}\in \mathcal{H}$ under $\hat{\mathbb{E}}$, if for each test function $\varphi\in C_{l,Lip}(\mathbb{R}^{m}\times\mathbb{R}^{n})$,
\begin{equation}\hat{\mathbb{E}}[\varphi(\bm{X},\bm{Y})]=\hat{\mathbb{E}}\left[\hat{\mathbb{E}}[\varphi(\bm{x},\bm{Y})]|_{\bm{x}=\bm{X}}\right]\label{I2}\end{equation}
whenever $\bar{\varphi}(\bm{x}):=\hat{\mathbb{E}}[|\varphi(\bm{x},\bm{Y})|]<\infty$ for all $\bm{x}$ and $\hat{\mathbb{E}}[|\bar{\varphi}(\bm{X})|]<\infty$.
A sequence $\{X_{n}, n\geq1\}$ of random variables is said to be independent if
\begin{eqnarray*}\hat{\mathbb{E}}[\varphi(X_{1},\cdots,X_{n},X_{n+1})]=\hat{\mathbb{E}}\left[\hat{\mathbb{E}}[\varphi(x_{1},\cdots,x_{n},X_{n+1})]|_{x_{1}=X_{1},\cdots,x_{n}=X_{n}}\right]\end{eqnarray*}
for all $n\geq1$ and $\varphi\in C_{l.Lip}(\mathbb{R}^{n+1})$.

Let $\bm{X_{1}}$ and $\bm{X_{2}}$ be two $n$-dimensional random vectors defined, respectively in sub-linear expectation spaces $(\Omega_{1},\mathcal{H}_{1},\hat{\mathbb{E}}_{1})$ and $(\Omega_{2},\mathcal{H}_{2},\hat{\mathbb{E}}_{2})$. They are called identically distributed, if
$$\hat{\mathbb{E}}_{1}[\varphi(\bm{X_{1}})]=\hat{\mathbb{E}}_{2}[\varphi(\bm{X_{2}})],\textrm{ for any }\varphi\in C_{l,Lip}(\mathbb{R}^{n}),$$
whenever the upper expectations are finite.} Under this framework, many extensions of dependence structures and limit theorems have
been investigated successively. We refer the readers to Peng (2008b, 2010) for the central limit theorem
and weak law of large numbers; Chen and Hu (2014), and Zhang (2016a) for the law of the iterated
algorithm; Zhang (2015) for the small derivation and Chung's
law of the iterated logarithm; Zhang (2016b) for the the moment inequalities of the maximum
partial sums and the Kolomogov type strong law of large numbers; Wu and Jiang (2018) for strong law of large numbers and Chover's law of the iterated logarithm; Wu et al. (2018) for the asymptotic approximation of inverse moment; Tang et al. (2019) for some exponential inequalities; Xi et al. (2019) for complete convergence and its application in nonparametric regression model; Zhang et al. (2020) for some results on convergence of the sums of independent random variables; Kuczmaszewska (2020) for complete convergence; Wu and Jiang (2020) for complete convergence and complete moment convergence of negatively dependent random variables, and so on. Since many basic properties or powerful tools such as various inequalities in classical probability space are no longer available in sub-linear expectation space, the study of limit theorems under sub-linear expectations is still a challenging work.

Zhang (2022) introduced the concept of extended negative dependence as follows.

\noindent\textbf{Definition 1.2.}~~{\it In a sub-linear expectation space $(\Omega,\mathcal{H},\hat{\mathbb{E}})$, a sequence $\{X_{n}, n\geq1\}$ of random variables is said to be upper (resp. lower) extended
negatively dependent if there exists some dominating constant $K\geq1$ such that
\begin{eqnarray*}\hat{\mathbb{E}}\left[\prod_{i=1}^{n}f_{i}(X_{i})\right]\leq K\prod_{i=1}^{n}\hat{\mathbb{E}}[f_{i}(X_{i})],\textrm{ for any }n\geq1,\end{eqnarray*}
whenever the non-negative functions $f_{i}\in C_{b,Lip}(\mathbb{R})$, $i=1,2,\ldots$, are all non-decreasing
(resp. all non-increasing). They are called extended negatively dependent (END, for short) if they are both
upper extended negatively dependent and lower extended negatively dependent.}

Next, we recall the capacities corresponding to the sub-linear expectations. Let $\mathcal{G}\subset\mathcal{F}$. A function $V: \mathcal{G}\rightarrow[0, 1]$ is called a capacity if
\begin{eqnarray*}V(\phi)=0,V(\Omega)=1,~and~V(A)\leq V(B),\textrm{ for any }A\subset B,~A,B\in\mathcal{G}.\end{eqnarray*}
It is said to be sub-additive if $V(A\bigcup B)\leq V(A)+V(B)$ for all $A,B\in \mathcal{G}$ with $A\bigcup B\in \mathcal{G}$.

In the sub-linear space $(\Omega,\mathcal{H},\hat{\mathbb{E}})$, a pair $(\mathbb{V},\mathcal{V})$ of capacities can be denoted by
\begin{eqnarray*}\mathbb{V}(A)=\inf\{\hat{\mathbb{E}}[\xi]:I_{A}\leq\xi,\xi\in\mathcal{H}\},\mathcal{V}(A)=1-\mathbb{V}(A^{c}),\textrm{ for }A\in\mathcal{F},\end{eqnarray*}
where $A^{c}$ is the complement set of $A$. Obviously, $\mathbb{V}$ is sub-additive, and
\begin{eqnarray*}\hat{\mathbb{E}}[f]\leq\mathbb{V}(A)\leq\hat{\mathbb{E}}[g],~\text{\Large $\hat{\varepsilon}$}[f]\leq\mathcal{V}(A)\leq\text{\Large $\hat{\varepsilon}$}[g],\textrm{ if }f\leq I_{A}\leq g\textrm{ and }f,g\in\mathcal{H}.\end{eqnarray*}
The inequalities above can conclude Markov's inequality in the same form of the classical probability space since $I(|X|\geq x)\leq |X|^{p}/x^{p}\in \mathcal{H}$. Also, the Choquet integrals/expecations $(C_{\mathbb{V}},C_{\mathcal{V}})$ can be defined by
\begin{eqnarray*}C_{V}[X]=\int_{0}^{\infty}V(X\geq t)dt+\int_{-\infty}^{0}[V(X\geq t)-1]dt\end{eqnarray*}
with $V$ being replaced by $\mathbb{V}$ and $\mathcal{V}$, respectively.

As stated above, the study under sub-linear expectations is quite challenging, one reason for which is the lack of powerful tools. For a sequence $\{X_{n},n\geq1\}$ of random variables where $X_{i}$ is negatively dependent to $(X_{i+1},\ldots,X_{n})$ for each $i=1,\ldots,n-1$, Zhang (2016b) proved that for any $1\leq p\leq2$,
$$\hat{\mathbb{E}}\left[\left|\max_{1\leq k\leq n}\sum_{i=1}^{k}X_{i}\right|^{p}\right]\leq2^{2-p}\sum_{i=1}^{n}\hat{\mathbb{E}}[|X_{i}|^{p}],$$
provided $\hat{\mathbb{E}}[X_n]\leq0$ for each $n\geq1$. However, as pointed out in Zhang (2022), both the independence in \eqref{I2} and the negative dependence raised in Zhang (2016b) have directions,
i.e., $Y$ is independent (resp. negatively dependent) relative to $X$ does not imply that $X$ is independent relative to $Y$ (resp. negatively dependent). The two concepts are restrictive and difficult to verify, that is why Zhang (2022) put forward the concepts of extended independence and extended negative dependence. Suppose that $\{X_{n},n\geq1\}$ is a sequence of END random variables under sub-linear expectations with $\hat{\mathbb{E}}[X_n]\leq0$ for each $n\geq1$, Zhang (2022) established the following inequality
\begin{equation}
\mathbb{V}\left(\sum_{i=1}^{n}X_i>x\right)\leq(1+Ke)x^{-2}\sum_{i=1}^{n}\hat{\mathbb{E}}[X_i^{2}].\label{I3}\end{equation}
It is not known whether the von Bahr-Esseen inequality \eqref{I1} holds for this dependence structure under sub-linear expectations since the concept of extended negative dependence is quite different from negative dependence. Our intention is that, maybe the exponent in \eqref{I3} can also be extended from 2 to $1<p<2$. Of course, it holds trivially for $0<p\leq1$. However, it is also a quite challenging work. First, the method used in Chen et al. (2014) does not work since the expression $E|X|=\int_{0}^{\infty}P(|X|>t)dt$ is unavailable for the sub-linear expectation $\hat{\mathbb{E}}[|X|]$. Second, if we adopt the method in Zhang (2022), the monotony of the function $F(x)=\frac{e^{x}-1-x}{|x|^{p}}$ can not be guaranteed for $1<p<2$. Therefore, in this work, we will proceed it with two steps. The von Bahr-Esseen type inequality is established for bounded END random variables by a function concerning the hyperbolic cosine borrowed from Hu et al. (2015), and then it is proved for unbounded END random variables under sub-linear expectations. By virtue of the inequality, we obtain the Kolmogorov type weak law of large numbers for partial sums of END random variables.  Another application concerns two results on complete convergence for weighted sums, which improve some existing ones in sub-linear expectation space.

Throughout this paper, let $C$ be a positive constant not depending on $n$, which may not be the same in various places. $a_{n}=O(b_{n})$ means $a_{n}\leq Cb_{n}$ for all $n\geq1$. Denote $x_{+}=xI(x>0)$, where $I(A)$ is the indicator function of the event $A$.

\section{Main results}
In this section, we first establish the von Bahr-Esseen type inequality for bounded END random variables under sub-linear expectations as follows.

\noindent\textbf {Theorem 2.1.}~~{\it Let $1<p\leq2$ and $\epsilon>0$. If $\{Z_n, n \geq 1\}$ is a sequence of upper END random variables in $(\Omega,\mathcal{H},\hat{\mathbb{E}})$ with $\hat{\mathbb{E}}[Z_n]\leq0$ and $|Z_n|\leq\epsilon$ for each $n\geq1$, then
\begin{eqnarray}
\mathbb{V}\left(\sum_{i=1}^{n}Z_i>\epsilon\right)\leq Ke\epsilon^{-p}\sum_{i=1}^{n}\hat{\mathbb{E}}[|Z_i|^{p}].\label{M1}
\end{eqnarray}
If $\{Z_n, n \geq 1\}$ is a sequence of END random variables in $(\Omega,\mathcal{H},\hat{\mathbb{E}})$ with $\hat{\mathbb{E}}[Z_n]=\text{\Large $\hat{\varepsilon}$}[Z_{n}]=0$ and $|Z_n|\leq\epsilon$ for each $n\geq1$, then
\begin{eqnarray}
\mathbb{V}\left(\left|\sum_{i=1}^{n}Z_i\right|>\epsilon\right)\leq 2Ke\epsilon^{-p}\sum_{i=1}^{n}\hat{\mathbb{E}}[|Z_i|^{p}].\label{M2}
\end{eqnarray}
}

Applying Theorem 2.1, the von Bahr-Esseen type inequality for unbounded END random variables under sub-linear expectations can be obtained as follows.

\noindent\textbf {Theorem 2.2.}~~{\it Let $1<p\leq2$. If $\{X_n, n \geq 1\}$ is a sequence of upper END random variables in $(\Omega,\mathcal{H},\hat{\mathbb{E}})$ with $\hat{\mathbb{E}}[X_n]\leq0$ for each $n\geq1$, then for any $x>0$,
\begin{eqnarray}
\mathbb{V}\left(\sum_{i=1}^{n}X_i>x\right)\leq 4^{p}(1+Ke)x^{-p}\sum_{i=1}^{n}\hat{\mathbb{E}}[|X_i|^{p}].\label{M3}
\end{eqnarray}
If $\{X_n, n \geq 1\}$ is a sequence of END random variables in $(\Omega,\mathcal{H},\hat{\mathbb{E}})$ with $\hat{\mathbb{E}}[X_n]=\text{\Large $\hat{\varepsilon}$}[X_{n}]=0$ for each $n\geq1$, then for any $x>0$,
\begin{eqnarray}
\mathbb{V}\left(\left|\sum_{i=1}^{n}X_i\right|>x\right)\leq 2^{2p+1}(1+Ke)x^{-p}\sum_{i=1}^{n}\hat{\mathbb{E}}[|X_i|^{p}].\label{M4}
\end{eqnarray}
}

By virtue of Theorem 2.2, we can obtain the following Kolmogorov type weak law of large numbers for END random variables under sub-linear expectations.

\noindent\textbf{Theorem 2.3.}~~{\it Let $\{X_{n},n\geq1\}$ be a sequence of END random variables in $(\Omega,\mathcal{H},\hat{\mathbb{E}})$. If $\lim_{c\rightarrow\infty}\sup_{i\geq1}\hat{\mathbb{E}}[(|X_{i}|-c)_{+}]=0$, then for any $\varepsilon>0$,
\begin{eqnarray}\mathbb{V}\left(\left\{\frac{\sum_{i=1}^{n}(X_{i}-\text{\Large $\hat{\varepsilon}$}[X_{i}])}{n}<-\varepsilon\right\}\bigcup\left\{\frac{\sum_{i=1}^{n}(X _{i}-\hat{\mathbb{E}}[X_{i}])}{n}>\varepsilon\right\}\right)\rightarrow0\textrm{ as }n\rightarrow\infty.\label{M5}\end{eqnarray}
Moreover, if $\hat{\mathbb{E}}[X_{n}]=\bar{\mu}$ and $\text{\Large $\hat{\varepsilon}$}[X_{n}]=\underline{\mu}$, then
\begin{eqnarray*}\mathcal{V}\left(\underline{\mu}-\varepsilon\leq\frac{\sum_{i=1}^{n}X_{i}}{n}\leq\bar{\mu}+\varepsilon\right)\rightarrow1\textrm{ as }n\rightarrow\infty.\end{eqnarray*}}

\noindent\textbf{Remark 2.1.}~~Zhang (2022) obtained the result where $\{X_{n},n\geq1\}$ is identically distributed. Hence, Theorem 2.3 extends the corresponding result of Zhang (2022) to the case of non-identical distribution.

Recall that a real valued function $l(x)$, positive and measurable on $\mathbb{R}^+$, is said to be slowly varying if for each $\lambda>0$,
$$\lim_{x\rightarrow\infty}\frac{l(\lambda x)}{l(x)}=1.$$

Another application of Theorem 2.2 is the following result on complete convergence for weighted sums of END random variables.

\noindent\textbf{Theorem 2.4.}~~{\it Let $0<p<2$, $\alpha>0$ and $\alpha p>1$. Let $\{X;X_{n},n\geq1\}$ be a sequence of identically distributed END random variables defined on $(\Omega,\mathcal{H},\hat{\mathbb{E}})$. Assume further that $\hat{\mathbb{E}}[X]=0$ and $\hat{\mathbb{E}}[|X|^{p}]<\infty$ if $p>1$. Let $l(x)$ be a slowly varying function and $\{a_{ni},1\leq i\leq n,n\geq1\}$ be an array of nonnegative real numbers satisfying $\sum_{i=1}^{n}a_{ni}^{q}=O(n)$ for some $q>p$. If $C_{\mathbb{V}}\left[|X|^{p}l(|X|^{1/\alpha})\right]<\infty$, then for any $\varepsilon>0$,
\begin{eqnarray}\sum_{n=1}^{\infty}n^{\alpha p-2}l(n)\mathbb{V}\left(\sum_{i=1}^{n}a_{ni}X_{i}>\varepsilon n^{\alpha}\right)<\infty.\label{M6}\end{eqnarray}
If we further assume $\text{\Large $\hat{\varepsilon}$}[X]=0$ when $p>1$, then for any $\varepsilon>0$,
\begin{eqnarray}\sum_{n=1}^{\infty}n^{\alpha p-2}l(n)\mathbb{V}\left(\left|\sum_{i=1}^{n}a_{ni}X_{i}\right|>\varepsilon n^{\alpha}\right)<\infty.\label{M7}\end{eqnarray}}

For $\alpha p=1$, we can obtain the following conclusion under some slightly different assumptions.

\noindent\textbf{Theorem 2.5.}~~{\it Let $0<p<2$ and $\{X;X_{n},n\geq1\}$ be a sequence of identically distributed END random variables defined on $(\Omega,\mathcal{H},\hat{\mathbb{E}})$. Assume further that $\hat{\mathbb{E}}[X]=0$ and $\lim_{c\rightarrow\infty}\hat{\mathbb{E}}[(|X|^{p}-c)_{+}]=0$ if $p\geq1$. Let $l(x)$ be a slowly varying function and $\{a_{ni},1\leq i\leq n,n\geq1\}$ be an array of nonnegative real numbers satisfying $\sum_{i=1}^{n}a_{ni}^{q}=O(n)$ for some $q>p$. If $C_{\mathbb{V}}\left[|X|^{p}l(|X|^{p})\right]<\infty$, then for any $\varepsilon>0$,
\begin{eqnarray}\sum_{n=1}^{\infty}n^{-1}l(n)\mathbb{V}\left(\sum_{i=1}^{n}a_{ni}X_{i}>\varepsilon n^{1/p}\right)<\infty.\label{M8}\end{eqnarray}
If we further assume $\text{\Large $\hat{\varepsilon}$}[X]=0$ when $p\geq1$, then for any $\varepsilon>0$,
\begin{eqnarray*}\sum_{n=1}^{\infty}n^{-1}l(n)\mathbb{V}\left(\left|\sum_{i=1}^{n}a_{ni}X_{i}\right|>\varepsilon n^{1/p}\right)<\infty.\end{eqnarray*}}

\noindent\textbf{Remark 2.2.}~~Zhong and Wu (2017) also obtained the corresponding results on complete convergence for weighted sums. It is easy to see that the assumption $\sum_{i=1}^{n}a_{ni}^{2}=O(n)$ is improved to $\sum_{i=1}^{n}a_{ni}^{q}=O(n)$ for some $q>p$. Moreover, they assumed that $l(x)$ is a monotone nondecreasing function which is not needed here. For $\alpha p>1$ and $1<p<2$, the hypothesis $\hat{\mathbb{E}}[|X|^{p}l(|X|^{1/\alpha})]\leq C_{\mathbb{V}}\left[|X|^{p}l(|X|^{1/\alpha})\right]$ in their result is stronger than $\hat{\mathbb{E}}[|X|^{p}]<\infty$ in Theorem 2.4. Therefore, our theorems improve the related results in Zhong and Wu (2017).
 
\section{Some lemmas}

In this section, we will present some lemmas which are essential in proving the results on complete convergence. The first lemma concerning the slowly varying function can be seen in Bai and Su (1985).

\noindent\textbf {Lemma 3.1.}~~If $l(x)$ is a slowly varying function, then

\noindent (i) $\lim_{k\rightarrow\infty}\sup_{2^{k}\leq x<2^{k+1}}\frac{l(x)}{l(2^k)}=1$;

\noindent (ii) for every $r>0$, $C_{1}2^{kr}l(2^{k})\leq \sum_{j=1}^{k}2^{jr}l(2^{j})\leq C_{2}2^{kr}l(2^{k})$;

\noindent (iii) for every $r<0$, $C_{1}2^{kr}l(2^{k})\leq \sum_{j=k}^{\infty}2^{jr}l(2^{j})\leq C_{2}2^{kr}l(2^{k})$.

\noindent The following lemma comes from Zhong and Wu (2017).

\noindent\textbf {Lemma 3.2.}~~Suppose that $X\in \mathcal{H}$, $p>0$, and $\alpha>0$, $l(x)$ is a slowly varying function.

\noindent (i) For any $c_0>0$,
$$C_{\mathbb{V}}\left[|X|^{p}l(|X|^{1/\alpha})\right]<\infty\Leftrightarrow\sum_{n=1}^{\infty}n^{\alpha p-1}l(n)\mathbb{V}(|X|>c_0n^{\alpha})<\infty.$$

\noindent (ii) If $C_{\mathbb{V}}\left[|X|^{p}l(|X|^{1/\alpha})\right]<\infty$, then for any $\theta>1$ and $c_0>0$,
$$\sum_{k=1}^{\infty}\theta^{k\alpha p}l(\theta^k)\mathbb{V}(|X|>c_0\theta^{k\alpha})<\infty.$$

Since the discontinuous function $I(x)\notin C_{l,Lip}(\mathbb{R})$, it should be transformed under the sub-linear expectations. Therefore, before presenting the last lemma, some notations are indispensable.

For $0<\mu<1$, let $g(x)\in C_{l,Lip}(\mathbb{R})$ be an even and non-increasing function on $[0,\infty)$ such that $0\leq g(x)\leq1$ for all $x$ and $g(x)=1$ if $|x|\leq\mu$, $g(x)=0$ if $|x|>1$. Hence,
\begin{eqnarray}I(|x|\leq\mu)\leq g(x)\leq I(|x|\leq1),~~I(|x|>1)\leq 1-g(x)\leq I(|x|>\mu).\label{L1}\end{eqnarray}

Let $g_{j}(x)\in C_{l,Lip}(\mathbb{R}),~j\geq1$ such that $0\leq g_{j}(x)\leq1$ for all $x$ and $g_{j}(\frac{x}{2^{j\alpha}})=1$ if $2^{(j-1)\alpha}<|x|\leq 2^{j\alpha}$, $g_{j}(\frac{x}{2^{j\alpha}})=0$ if $|x|\leq \mu2^{(j-1)\alpha}$ or $|x|>(1+\mu)2^{j\alpha}$. Thus,
\begin{eqnarray}g_{j}\left(\frac{x}{2^{j\alpha}}\right)\leq I\left(\mu2^{(j-1)\alpha}<|x|\leq(1+\mu)2^{j\alpha}\right),\label{L2}\end{eqnarray}
and
\begin{eqnarray}|x|^{t}g\left(\frac{x}{2^{k\alpha}}\right)\leq 1+\sum_{j=1}^{k}|x|^{t}g_{j}\left(\frac{x}{2^{j\alpha}}\right),\textrm{ for any }t>0.\label{L3}\end{eqnarray}

\noindent\textbf {Lemma 3.3.}~~Suppose that $X\in \mathcal{H}$, $0<\mu<1$, $\alpha>0$ and $p>0$. If $C_{\mathbb{V}}\left[|X|^{p}l(|X|^{1/\alpha})\right]<\infty$, then for any $s>p$,
$$\sum_{n=1}^{\infty}n^{\alpha p-\alpha s-1}l(n)\hat{\mathbb{E}}\left[|X|^{s}g\left(\frac{\mu X}{n^{\alpha}}\right)\right]<\infty.$$
\noindent\textbf{Proof.}~~It follows from \eqref{L2}, \eqref{L3}, Lemmas 3.1 and 3.2 that
\begin{eqnarray*}&&\sum_{k=0}^{\infty}\sum_{2^{k}\leq n<2^{k+1}}n^{\alpha p-\alpha s-1}l(n)\hat{\mathbb{E}}\left[|X|^{s}g\left(\frac{\mu X}{n^{\alpha}}\right)\right]\\
&\leq&\sum_{k=0}^{\infty}2^{\alpha(p-s)k}l(2^{k})\hat{\mathbb{E}}\left[|X|^{s}g\left(\frac{\mu X}{2^{\alpha(k+1)}}\right)\right]\\
&\leq&C\sum_{k=0}^{\infty}2^{\alpha(p-s)k}l(2^{k})\hat{\mathbb{E}}\left[1+\sum_{j=0}^{k}|X|^{s}g_{j+1}\left(\frac{\mu X}{2^{\alpha(j+1)}}\right)\right]\\
&\leq&C\sum_{k=0}^{\infty}2^{\alpha(p-s)k}l(2^{k})+C\sum_{k=0}^{\infty}2^{\alpha(p-s)k}l(2^{k})\sum_{j=0}^{k}\hat{\mathbb{E}}\left[|X|^{s}g_{j+1}\left(\frac{\mu X}{2^{\alpha(j+1)}}\right)\right]\\
&\leq&C+C\sum_{j=0}^{\infty}\hat{\mathbb{E}}\left[|X|^{s}g_{j+1}\left(\frac{\mu X}{2^{\alpha(j+1)}}\right)\right]\sum_{k=j}^{\infty}2^{\alpha(p-s)k}l(2^{k})\\
&\leq&C+C\sum_{j=0}^{\infty}2^{\alpha pj}l(2^{j})\mathbb{V}(|X|>\mu^{2}2^{j\alpha})<\infty.
\end{eqnarray*}
The proof is completed. $\Box$

\section{The proofs of the main results}

\noindent\textbf{Proof of Theorem 2.1.}~~Note that for any $t>0$, the non-negative function $f(x)=e^{tx}\in C_{b,Lip}(\mathbb{R})$ since $f'(x)=te^{tx}\leq te^{t\epsilon}$ for any $|x|\leq\epsilon$. It follows from Markov's inequality and Definition 1.2 that
\begin{eqnarray}\mathbb{V}\left(\sum_{i=1}^{n}Z_{i}>\epsilon\right)\leq e^{-t\epsilon}\hat{\mathbb{E}}\left[\exp\left(t\sum_{i=1}^{n}Z_{i}\right)\right]=e^{-t\epsilon}\hat{\mathbb{E}}\left[\prod_{i=1}^{n}e^{tZ_{i}}\right]\leq Ke^{-tx}\prod_{i=1}^{n}\hat{\mathbb{E}}[e^{tZ_{i}}].\label{P1}\end{eqnarray}
Define the function $H(x): \mathbb{R}\rightarrow\mathbb{R}$ by
\begin{equation*}
H(x)=\begin{cases}
     {{\cosh x-1}\over{|x|^{p}}},  & if~x\neq 0, \\
     {{1}\over{2}},  & if~x=0~and~p=2, \\
     0,  & if~x=0~and~0<p<2.
\end{cases}
\end{equation*}
It is easy to verify that $H(x)$ is a nonnegative, continuous, even function, and is increasing when $x\geq0$. Thus, for all $|Z_{i}|\leq\epsilon$, we have $H(tZ_{i})\leq H(t\epsilon)$, $1\leq i\leq n$. Therefore, we obtain by $e^x-1-x\leq 2(\cosh x-1)$ for any $x\in\mathbb{R}$ that
\begin{eqnarray}e^{tZ_{i}}&=&1+tZ_{i}+\frac{e^{tZ_{i}}-1-tZ_{i}}{|tZ_{i}|^p}|tZ_{i}|^p\nonumber\\
&\leq&1+tZ_{i}+\frac{2\big(\cosh (tZ_{i})-1\big)}{|tZ_{i}|^p}|tZ_{i}|^p\nonumber\\
&\leq&1+tZ_{i}+\frac{2\big(\cosh (t\epsilon)-1\big)}{\epsilon^p}|Z_{i}|^p.\label{P2}\end{eqnarray}
Note that $\hat{\mathbb{E}}[Z_{i}]\leq0$ and $1+x\leq e^x$ for any $x\in\mathbb{R}$. A combination of \eqref{P1} and \eqref{P2} gives that
\begin{eqnarray}\mathbb{V}\left(\sum_{i=1}^{n}Z_{i}>\epsilon\right)&\leq&Ke^{-t\epsilon}\prod_{i=1}^{n}\left[1+t\hat{\mathbb{E}}[Z_{i}]+\frac{2\big(\cosh (t\epsilon)-1\big)}{\epsilon^p}\hat{\mathbb{E}}[|Z_{i}|^p]\right]\nonumber\\
&\leq&Ke^{-t\epsilon}\prod_{i=1}^{n}\left[1+\frac{2\big(\cosh (t\epsilon)-1\big)}{\epsilon^p}\hat{\mathbb{E}}[|Z_{i}|^p]\right]\nonumber\\
&\leq&Ke^{-t\epsilon}\prod_{i=1}^{n}\exp\left\{\frac{2\big(\cosh (t\epsilon)-1\big)}{\epsilon^p}\hat{\mathbb{E}}[|Z_{i}|^p]\right\}\nonumber\\
&=&K\exp\left\{-t\epsilon+\frac{2\big(\cosh (t\epsilon)-1\big)}{\epsilon^p}\sum_{i=1}^{n}\hat{\mathbb{E}}[|Z_{i}|^p]\right\}.\label{P3}\end{eqnarray}
Since $t>0$ is arbitrary, we take $t=\epsilon^{-1}\ln\big(1+\frac{\epsilon^{p}}{\sum_{i=1}^{n}\hat{\mathbb{E}}[|Z_{i}|^p]}\big)$ in \eqref{P3} to obtain that
\begin{eqnarray*}\mathbb{V}\left(\sum_{i=1}^{n}Z_{i}>\epsilon\right)&\leq&K\exp\left\{-\ln\left(1+\frac{\epsilon^{p}}{\sum_{i=1}^{n}\hat{\mathbb{E}}[|Z_{i}|^p]}\right)+\frac{\epsilon^{p}}{\epsilon^{p}+\sum_{i=1}^{n}\hat{\mathbb{E}}[|Z_{i}|^p]}\right\}\\
&\leq&Ke\frac{\sum_{i=1}^{n}\hat{\mathbb{E}}[|Z_{i}|^p]}{\epsilon^{p}+\sum_{i=1}^{n}\hat{\mathbb{E}}[|Z_{i}|^p]}\leq Ke\epsilon^{-p}\sum_{i=1}^{n}\hat{\mathbb{E}}[|Z_{i}|^p],\end{eqnarray*}
which proves \eqref{M1}. Replacing $Z_{i}$ with $-Z_{i}$ for each $1\leq i\leq n$ in \eqref{M1} and noting that
$$\mathbb{V}\left(\left|\sum_{i=1}^{n}Z_i\right|>\epsilon\right)\leq\mathbb{V}\left(\sum_{i=1}^{n}Z_i>\epsilon\right)+\mathbb{V}\left(\sum_{i=1}^{n}(-Z_i)>\epsilon\right),$$
we get \eqref{M2} immediately. The details are omitted. $\Box$

\noindent\textbf{Proof of Theorem 2.2.}~~Since the desired result holds trivially if $4^{p}(1+Ke)x^{-p}\sum_{i=1}^{n}\hat{\mathbb{E}}[|X_i|^{p}]>1$, we only need to prove it under the assumption $4^{p}(1+Ke)x^{-p}\sum_{i=1}^{n}\hat{\mathbb{E}}[|X_i|^{p}]\leq1$ in what follows. Define for each $1\leq i\leq n$ and any $x>0$ that
\begin{eqnarray*}X_{i}(1)&=&-\frac{x}{4}I\big(X_{i}<-\frac{x}{4}\big)+X_{i}I\big(|X_{i}|\leq \frac{x}{4}\big)+\frac{x}{4}I\big(X_{i}>\frac{x}{4}\big);\\
X_{i}(2)&=&\big(X_{i}+\frac{x}{4}\big)I\big(X_{i}<-\frac{x}{4}\big)+\big(X_{i}-\frac{x}{4}\big)I\big(X_{i}>\frac{x}{4}\big).\end{eqnarray*}
It is easy to see that $X_{i}=X_{i}(1)+X_{i}(2)$ for each $1\leq i\leq n$ and $\{X_{i}(1)-EX_{i}(1),i\geq1\}$ is still a sequence of END random variables since $f_{\frac{x}{4}}(x)=\max\{-\frac{x}{4},\min(x,\frac{x}{4})\}\in C_{l,Lip}(\mathbb{R})$. Note that $|X_{i}(1)-\hat{\mathbb{E}}[X_{i}(1)]|\leq\frac{x}{2}$ and $\hat{\mathbb{E}}[X_{i}(1)-\hat{\mathbb{E}}[X_{i}(1)]]=0$. It follows from Theorem 2.1, $|X_{i}(1)|\leq|X_{i}|$, $C_{r}$-inequality and Jensen's inequality that
\begin{eqnarray}\mathbb{V}\left(\sum_{i=1}^{n}(X_{i}(1)-\hat{\mathbb{E}}[X_{i}(1)])>\frac{x}{2}\right)&\leq&Ke\big(\frac{x}{2}\big)^{-p}\sum_{i=1}^{n}\hat{\mathbb{E}}[|X_{i}(1)-\hat{\mathbb{E}}[X_{i}(1)]|^p]\nonumber\\
&\leq&Ke\big(\frac{x}{2}\big)^{-p}2^{p-1}\sum_{i=1}^{n}\big(\hat{\mathbb{E}}[|X_{i}|^{p}]+(\hat{\mathbb{E}}[|X_{i}|])^p\big)\nonumber\\
&\leq&Ke\big(\frac{x}{2}\big)^{-p}2^{p}\sum_{i=1}^{n}\hat{\mathbb{E}}[|X_{i}|^{p}]\nonumber\\
&=&4^{p}Kex^{-p}\sum_{i=1}^{n}\hat{\mathbb{E}}[|X_{i}|^{p}].\label{P4}\end{eqnarray}
By Markov's inequality we can easily obtain that
\begin{eqnarray}\mathbb{V}\left(\sum_{i=1}^{n}X_{i}(2)>\frac{x}{4}\right)\leq\mathbb{V}\left(\bigcup_{i=1}^{n}\Big\{|X_{i}|>\frac{x}{4}\Big\}\right)\leq\sum_{i=1}^{n}\mathbb{V}\left(|X_{i}|>\frac{x}{4}\right)
\leq4^{p}x^{-p}\sum_{i=1}^{n}\hat{\mathbb{E}}[|X_{i}|^{p}].\label{P5}\end{eqnarray}
Furthermore, observing that $(|X_{i}|-\frac{x}{4})I(|X_{i}|>\frac{x}{4})\leq|X_{i}|I(|X_{i}|>\frac{x}{4})\leq(\frac{4}{x})^{p-1}|X_{i}|^{p}$ and $4^{p}(1+Ke)x^{-p}\sum_{i=1}^{n}\hat{\mathbb{E}}[|X_i|^{p}]\leq1$, we have
\begin{eqnarray}\sum_{i=1}^{n}\big(\hat{\mathbb{E}}[X_{i}(1)]-\hat{\mathbb{E}}[X_{i}]\big)&\leq&\sum_{i=1}^{n}\hat{\mathbb{E}}\big[|X_{i}(1)-X_{i}|\big]\nonumber\\
&\leq&\sum_{i=1}^{n}\hat{\mathbb{E}}\Big[\big(|X_{i}|-\frac{x}{4}\big)I\big(|X_{i}|>\frac{x}{4}\big)\Big]\nonumber\\
&\leq&\big(\frac{4}{x}\big)^{p-1}\sum_{i=1}^{n}\hat{\mathbb{E}}[|X_{i}|^{p}]\nonumber\\
&=&4^{p}x^{-p}\sum_{i=1}^{n}\hat{\mathbb{E}}[|X_{i}|^{p}]\cdot\frac{x}{4}\leq\frac{x}{4}.\label{P6}\end{eqnarray}
Recall that $\hat{\mathbb{E}}[X_i]\leq0$ for each $1\leq i\leq n$. It follows from \eqref{P4}-\eqref{P6} that
\begin{eqnarray*}\mathbb{V}\left(\sum_{i=1}^{n}X_{i}>x\right)&\leq&
\mathbb{V}\left(\sum_{i=1}^{n}(X_{i}(1)-\hat{\mathbb{E}}[X_{i}(1)])+\sum_{i=1}^{n}\big(\hat{\mathbb{E}}[X_{i}(1)]-\hat{\mathbb{E}}[X_{i}]\big)+\sum_{i=1}^{n}X_{i}(2)>x\right)\\
&\leq&\mathbb{V}\left(\sum_{i=1}^{n}(X_{i}(1)-\hat{\mathbb{E}}[X_{i}(1)])>\frac{x}{2}\right)+\mathbb{V}\left(\sum_{i=1}^{n}X_{i}(2)>\frac{x}{4}\right)\\
&\leq&4^{p}(1+Ke)x^{-p}\sum_{i=1}^{n}\hat{\mathbb{E}}[|X_i|^{p}],\end{eqnarray*}
which implies \eqref{M3}. Replacing $X_{i}$ with $-X_{i}$ for each $1\leq i\leq n$ in \eqref{M3}, we can easily obtain \eqref{M4}. The proof is completed. $\Box$

\noindent\textbf{Proof of Theorem 2.3.}~~Define for each $1\leq i\leq n$ that
\begin{eqnarray*}X_{i}(3)&=&-cI(X_{i}<-c)+X_{i}I(|X_{i}|\leq c)+cI(X_{i}>c),\\
X_{i}(4)&=&(X_{i}+c)I(X_{i}<-c)+(X_{i}-c)I(X_{i}>c),\end{eqnarray*}
where $c>0$ is sufficiently large such that $\sup_{i\geq1}\hat{\mathbb{E}}[(|X_{i}|-c)_{+}]\leq\varepsilon/3$. It is easy to see that $X_{i}=X_{i}(3)+X_{i}(4)$ for each $1\leq i\leq n$ and $\{X_{i}(3)-\hat{\mathbb{E}}[X_{i}(3)],i\geq1\}$ is a sequence of END random variables. Moreover, we can also easily check that
\begin{eqnarray*}n^{-1}\sum_{i=1}^{n}\big(\hat{\mathbb{E}}[X_{i}(3)]-\hat{\mathbb{E}}[X_{i}]\big)&\leq&n^{-1}\sum_{i=1}^{n}\hat{\mathbb{E}}\big[|X_{i}(3)-X_{i}|\big]\nonumber\\
&\leq&n^{-1}\sum_{i=1}^{n}\hat{\mathbb{E}}\big[(|X_{i}|-c)_{+}\big]\nonumber\\
&\leq&\sup_{i\geq1}\hat{\mathbb{E}}[(|X_{i}|-c)_{+}]\leq\varepsilon/3.\end{eqnarray*}
Hence, we obtain by Theorem 2.2 and Markov's inequality that
\begin{eqnarray*}&&\mathbb{V}\left(\sum_{i=1}^{n}(X _{i}-\hat{\mathbb{E}}[X_{i}])>\varepsilon n\right)\\
&=&
\mathbb{V}\left(\sum_{i=1}^{n}(X_{i}(3)-\hat{\mathbb{E}}[X_{i}(3)])+\sum_{i=1}^{n}\big(\hat{\mathbb{E}}[X_{i}(3)]-\hat{\mathbb{E}}[X_{i}]\big)+\sum_{i=1}^{n}X_{i}(4)>\varepsilon n\right)\\
&\leq&\mathbb{V}\left(\sum_{i=1}^{n}(X_{i}(3)-\hat{\mathbb{E}}[X_{i}(3)])>\frac{\varepsilon n}{3}\right)+\mathbb{V}\left(\sum_{i=1}^{n}X_{i}(4)>\frac{\varepsilon n}{3}\right)\\
&\leq&\big(\frac{3}{\varepsilon}\big)^{2}n^{-2}\sum_{i=1}^{n}\hat{\mathbb{E}}\big[(X_{i}(3)-\hat{\mathbb{E}}[X_{i}(3)])^{2}\big]+\mathbb{V}\left(\sum_{i=1}^{n}|X_{i}(4)|>\varepsilon n/3\right)\\
&\leq&\big(\frac{3}{\varepsilon}\big)^{2}n^{-2}\sum_{i=1}^{n}\hat{\mathbb{E}}\big[(X_{i}(3))^{2}\big]+\frac{3}{\varepsilon}n^{-1}\sum_{i=1}^{n}\hat{\mathbb{E}}[|X_{i}(4)|]\\
&\leq&\big(\frac{3}{\varepsilon}\big)^{2}n^{-1}c^{2}+\frac{3}{\varepsilon}\sup_{i\geq1}\hat{\mathbb{E}}[(|X_{i}|-c)_{+}]\rightarrow0\end{eqnarray*}
by letting $n\rightarrow\infty$ first and then $c\rightarrow\infty$. Hence, the desired result \eqref{M5} has been proved. $\Box$

\noindent\textbf{Proof of Theorem 2.4.}~~It follows from H\"{o}lder's inequality that $\sum_{i=1}^{n}a_{ni}^{\gamma}=O(n)$ for any $0<\gamma<q$. The proof will be conducted under the following two cases.

\textbf{Case 1.} $0<p<1$.

It is easy to check that
\begin{eqnarray*}&&\sum_{n=1}^{\infty}n^{\alpha p-2}l(n)\mathbb{V}\left(\sum_{i=1}^{n}a_{ni}X_{i}>\varepsilon n^{\alpha}\right)\\
&\leq&\sum_{n=1}^{\infty}n^{\alpha p-2}l(n)\mathbb{V}\left(\sum_{i=1}^{n}a_{ni}X_{i}I(|X_{i}|\leq n^{\alpha})>\varepsilon n^{\alpha}\right)\\
&&+\sum_{n=1}^{\infty}n^{\alpha p-2}l(n)\mathbb{V}\left(\bigcup_{i=1}^{n}\{|X_{i}|>n^{\alpha}\}\right)\\
&=:&I_{1}+I_{2}.\end{eqnarray*}
Denote $q'=\min\{1,q\}$. Note by (\ref{L1}) that $I(|X_{i}|\leq n^{\alpha})\leq g\left(\frac{\mu X_{i}}{n^{\alpha}}\right)$. Hence, it follows from Markov's inequality and Lemma 3.3 that
\begin{eqnarray*}I_{1}&\leq&C\sum_{n=1}^{\infty}n^{\alpha p-\alpha q'-2}l(n)\sum_{i=1}^{n}a_{ni}^{q'}\hat{\mathbb{E}}\left[|X_{i}|^{q'}g\left(\frac{\mu X_{i}}{n^{\alpha}}\right)\right]\\
&\leq&C\sum_{n=1}^{\infty}n^{\alpha p-\alpha q'-1}l(n)\hat{\mathbb{E}}\left[|X|^{q'}g\left(\frac{\mu X}{n^{\alpha}}\right)\right]<\infty.\end{eqnarray*}
On the other hand, by \eqref{L1} and Lemma 3.2 we can obtain that
\begin{eqnarray*}I_{2}&\leq&\sum_{n=1}^{\infty}n^{\alpha p-2}l(n)\sum_{i=1}^{n}\mathbb{V}\left(|X_{i}|>n^{\alpha}\right)\\
&\leq&\sum_{n=1}^{\infty}n^{\alpha p-2}l(n)\sum_{i=1}^{n}\hat{\mathbb{E}}\left[1-g\left(\frac{X_{i}}{n^{\alpha}}\right)\right]\\
&=&\sum_{n=1}^{\infty}n^{\alpha p-1}l(n)\hat{\mathbb{E}}\left[1-g\left(\frac{X}{n^{\alpha}}\right)\right]\\
&\leq&\sum_{n=1}^{\infty}n^{\alpha p-1}l(n)\mathbb{V}\left(|X|>\mu n^{\alpha}\right)<\infty.\end{eqnarray*}

\textbf{Case 2.} $1\leq p<2$.

Define for each $1\leq i\leq n$, $n\geq1$ that
\begin{eqnarray*}X_{ni}=-n^{\alpha}I(X_{i}<-n^{\alpha})+X_{i}I(|X_{i}|\leq n^{\alpha})+n^{\alpha}I(X_{i}>n^{\alpha}).\end{eqnarray*}
It is easy to see that $\{X_{ni},1\leq i\leq n,n\geq1\}$ is still an array of END random variables under sub-linear expectations. We will firstly show
\begin{eqnarray}n^{-\alpha}\left|\sum_{i=1}^{n}a_{ni}\hat{\mathbb{E}}[X_{ni}]\right|\rightarrow0\textrm{ as }n\rightarrow\infty.\label{P7}\end{eqnarray}
Noting that $\hat{\mathbb{E}}[X]=0$ if $p>1$, and
$$(|X|-n^{\alpha})_{+}\leq|X|I(|X|>n^{\alpha})\leq n^{\alpha(1-p)}|X|^{p}I(|X|>n^{\alpha})\leq n^{\alpha(1-p)}|X|^{p},$$
we derive by $\alpha p>1$ and $\hat{\mathbb{E}}[|X|^{p}]<\infty$ that
\begin{eqnarray*}n^{-\alpha}\left|\sum_{i=1}^{n}a_{ni}\hat{\mathbb{E}}[X_{ni}]\right|&=&n^{-\alpha}\left|\sum_{i=1}^{n}a_{ni}(\hat{\mathbb{E}}[X_{ni}]-\hat{\mathbb{E}}[X_{i}])\right|\\
&\leq&n^{-\alpha}\sum_{i=1}^{n}a_{ni}\hat{\mathbb{E}}[|X_{i}-X_{ni}|]\\
&\leq&Cn^{1-\alpha}\hat{\mathbb{E}}[(|X|-n^{\alpha})_{+}]\\
&\leq&Cn^{1-\alpha p}\hat{\mathbb{E}}[|X|^{p}]\rightarrow0\textrm{ as }n\rightarrow\infty.\end{eqnarray*}
If $p=1$, we also have by $\alpha>1$ that
\begin{eqnarray*}n^{-\alpha}\left|\sum_{i=1}^{n}a_{ni}\hat{\mathbb{E}}[X_{ni}]\right|\leq n^{-\alpha}\sum_{i=1}^{n}a_{ni}\hat{\mathbb{E}}[|X_{ni}|]
\leq Cn^{1-\alpha}\hat{\mathbb{E}}[|X|]\rightarrow0\textrm{ as }n\rightarrow\infty.\end{eqnarray*}
Hence, \eqref{P7} has been proved, from which we can conclude that
\begin{eqnarray*}&&\sum_{n=1}^{\infty}n^{\alpha p-2}l(n)\mathbb{V}\left(\sum_{i=1}^{n}a_{ni}X_{i}>\varepsilon n^{\alpha}\right)\\
&\leq&C\sum_{n=1}^{\infty}n^{\alpha p-2}l(n)\mathbb{V}\left(\sum_{i=1}^{n}a_{ni}\big(X_{ni}-\hat{\mathbb{E}}[X_{ni}]\big)>\frac{\varepsilon}{2} n^{\alpha}\right)\\
&&+\sum_{n=1}^{\infty}n^{\alpha p-2}l(n)\mathbb{V}\left(\bigcup_{i=1}^{n}\{|X_{i}|>n^{\alpha}\}\right)\\
&=:&I_{3}+I_{2}.\end{eqnarray*}
Since $I_{2}<\infty$ has been proved in Case 1, we only need to show $I_{3}<\infty$. Take $\beta=\min\{2,q\}$. Note by (\ref{L1}) that
$$|X_{ni}|^{\beta}=|X_{i}|^{\beta}I(|X_{i}|\leq n^{\alpha})+n^{\alpha\beta}I(|X_{i}|>n^{\alpha})\leq|X_{i}|^{\beta}g\left(\frac{\mu X_{i}}{n^{\alpha}}\right)+n^{\alpha\beta}\left[1-g\left(\frac{X_{i}}{n^{\alpha}}\right)\right].$$
Therefore, we obtain by \eqref{L1} again, Theorem 2.2, $C_{r}$-inequality, Jensen's inequality, Lemma 3.2 and Lemma 3.3 that
\begin{eqnarray*}I_{3}&\leq&C\sum_{n=1}^{\infty}n^{\alpha p-\alpha\beta-2}l(n)\sum_{i=1}^{n}a_{ni}^{\beta}\hat{\mathbb{E}}\big[|X_{ni}-\hat{\mathbb{E}}[X_{ni}]|^{\beta}\big]\\
&\leq&C\sum_{n=1}^{\infty}n^{\alpha p-\alpha\beta-2}l(n)\sum_{i=1}^{n}a_{ni}^{\beta}\hat{\mathbb{E}}\big[|X_{ni}|^{\beta}\big]\\
&\leq&C\sum_{n=1}^{\infty}n^{\alpha p-\alpha\beta-2}l(n)\sum_{i=1}^{n}a_{ni}^{\beta}\left\{\hat{\mathbb{E}}\left[|X_{i}|^{\beta}g\left(\frac{\mu X_{i}}{n^{\alpha}}\right)\right]+n^{\alpha\beta}\hat{\mathbb{E}}\left[1-g\left(\frac{X_{i}}{n^{\alpha}}\right)\right]\right\}\\
&\leq&C\sum_{n=1}^{\infty}n^{\alpha p-\alpha\beta-1}l(n)\hat{\mathbb{E}}\left[|X|^{\beta}g\left(\frac{\mu X}{n^{\alpha}}\right)\right]+C\sum_{n=1}^{\infty}n^{\alpha p-1}l(n)\hat{\mathbb{E}}\left[1-g\left(\frac{X}{n^{\alpha}}\right)\right]\\
&\leq&C\sum_{n=1}^{\infty}n^{\alpha p-\alpha\beta-1}l(n)\hat{\mathbb{E}}\left[|X|^{\beta}g\left(\frac{\mu X}{n^{\alpha}}\right)\right]+C\sum_{n=1}^{\infty}n^{\alpha p-1}l(n)\mathbb{V}\left(|X|>\mu n^{\alpha}\right)\\
&<&\infty.\end{eqnarray*}
Therefore, (\ref{M6}) is proved. Replacing $X_{i}$ with $-X_{i}$ for each $1\leq i\leq n$ in \eqref{M6}, we can obtain \eqref{M7} without much effort. The proof is hence completed. $\Box$

\noindent\textbf{Proof of Theorem 2.5.}~~We only give the proof of \eqref{M8}. In view of the proof of Theorem 2.4, we only need to verify that \eqref{P7} holds under the assumption $\lim_{c\rightarrow\infty}\hat{\mathbb{E}}[(|X|^{p}-c)_{+}]=0$ when $1\leq p<2$ with $\alpha p=1$. Actually, it is easy to check that
\begin{eqnarray*}(|X|-n^{1/p})_{+}&\leq&|X|I(|X|>n^{1/p})\\
&\leq&n^{(1-p)/p}|X|^{p}I(|X|^{p}>n)\\
&\leq&2n^{(1-p)/p}(|X|^{p}-n/2)I(|X|^{p}>n)\\
&\leq&2n^{(1-p)/p}(|X|^{p}-n/2)_{+}.\end{eqnarray*}
Hence, we obtain by the assumption $\hat{\mathbb{E}}[X]=0$ for $p\geq1$ that
\begin{eqnarray*}n^{-1/p}\left|\sum_{i=1}^{n}a_{ni}\hat{\mathbb{E}}[X_{ni}]\right|&\leq&n^{-1/p}\sum_{i=1}^{n}a_{ni}\hat{\mathbb{E}}[|X_{i}-X_{ni}|]\\
&\leq&Cn^{1-1/p}\hat{\mathbb{E}}[(|X|-n^{1/p})_{+}]\\
&\leq&C\hat{\mathbb{E}}[(|X|^{p}-n/2)_{+}]\rightarrow0\textrm{ as }n\rightarrow\infty.\end{eqnarray*}
The proof is completed. $\Box$


\begin{thebibliography}{1}

\bibitem{RefJ1} Artzner Ph., Delbaen F., Eber J.M., Heath D., 1997. Thinking coherently. RISK, 10, 68-71.

\bibitem{RefJ2} Asadian N., Fakoor V., Bozorgnia A., 2006. Rosenthal's type inequalities for negatively orthant dependent random variables. Journal of the Iranian Statistical Society, 5, 69-75.

\bibitem{RefJ3} Bai Z.D., Su C., 1985. The complete convergence for partial sums of i.i.d. random variables. Science in China, Series A, 12, 31-47.

\bibitem{RefJ4} Bryc W., Smolenski W., 1993. Moment conditions for almost sure convergence of weakly correlated random variables. Proceedings of the American Mathematical Society, 119(2), 629-635.

\bibitem{RefJ5} Chatterji S.D., 1969. An $L^p$-convergence theorem. The Annals of Mathematical Statistics, 40, 1068-1070.

\bibitem{RefJ6} Chen P.Y., Bai P, Sung S.H., 2014. The von Bahr-Esseen moment inequality for pairwise independent random variables and applications. Journal of Mathematical Analysis and Applications, 419(2), 1290-1302.

\bibitem{RefJ7} Chen Z.J., Epstein L., 2002. Ambiguity, risk and asset returns in continuous time. Econometrica, 70(4), 1403-1443.

\bibitem{RefJ8} Chen Z.J., Hu F., 2014. A law of the iterated logarithm under sublinear expectations. Journal of Financial Engineering, 1(2), Article ID 1450015, 23 pages.

\bibitem{RefJ9} Denis L., Martini C., 2006. A theoretical framework for the pricing of contingent claims in the presence of model uncertainty. The Annals of Applied Probability, 16(2), 827-852.

\bibitem{RefJ10} El Karoui N., Peng S.G., Quenez M.C., 1997. Backward stochastic differential equation in finance. Mathematical Finance, 7(1), 1-71.

\bibitem{RefJ11} Follmer H., Schied A., 2004. Statistic Finance, An Introduction in Discrete Time (2nd Edition), Walter de Gruyter.

\bibitem{RefJ12} Gilboa I., 1987. Expected utility theory with purely subjective non-additive probabilities. Journal of Mathematical Economics, 16(1), 65-88.

\bibitem{RefJ13} Hu T.-C., Wang K.L., Rosalsky A., 2015. Complete convergence theorems for extended negatively dependent random variables. Sankhy\={a} A: The Indian Journal of Statistics, 77(1), 1-29.

\bibitem{RefJ14} Kuczmaszewska A., 2020. Complete convergence for widely acceptable random variables under the sublinear expectations. Journal of Mathematical Analysis and Applications, 484(1), Article ID 123662.

\bibitem{RefJ15} Marinacci M. 1999. Limit laws for non-additive probabilities and their frequentist interpretation. Journal of Economic Theory, 84(2), 145-195.

\bibitem{RefJ16} Peng S.G., 1997. Backward SDE and related $g$-expectations, in Backward Stochastic Differential Equations, Pitman Research Notes in Math. Series, No. 364, El Karoui Mazliak edit. 141-159.

\bibitem{RefJ17} Peng S.G.,1999. Monotonic limit theorem of BSDE and nonlinear decomposition theorem of Doob-Meyer type. Probability Theory and Related Fields, 113(4), 473-499.

\bibitem{RefJ18} Peng S.G., 2006. $G$-expectation, $G$-Brownian motion and related stochastic calculus of Itos type. The Abel Symposium 2005, Abel Symposia 2, Edit. Benth et. al., 541-567, Springer-Verlag, 2006.

\bibitem{RefJ19} Peng S.G., 2008a. Multi-dimensional $G$-Brownian motion and related stochastic calculus under $G$-expectation. Stochastic Processes and Their Applications, 118(12), 2223-2253.

\bibitem{RefJ20} Peng S.G., 2008b. A new central limit theorem under sublinear expectations. arXiv preprint arXiv:0803.2656.

\bibitem{RefJ21} Peng S.G., 2010. Nonlinear expectations and stochastic calculus under uncertainty. arXiv preprint arXiv:1002.4546.

\bibitem{RefJ22} Shao Q.M., 2000. A comparison theorem on moment inequalities between negatively associated and independent random variables. Journal of Theoretical Probability, 13(2), 343-356.

\bibitem{RefJ23} Tang X.F., Wang X.J., Wu Y., 2019. Exponential inequalities under sub-linear expectations with applications to strong law of large numbers. Filomat, 33(10), 2951-2961.

\bibitem{RefJ24} von Bahr B., Esseen C.G., 1965. Inequalities for the $r$th absolute moment of a sum of random variables, $1\leq r\leq2$. The Annals of Mathematical Statistics, 36(1), 299-303.

\bibitem{RefJ25} Wang X.J., Xu C., Hu T.-C., Volodin A., Hu S.H., 2014. On complete convergence for widely orthant dependent random variables and its applications in nonparametric regression models. TEST, 23,
607-629.

\bibitem{RefJ26} Wu Q.Y., Jiang Y.Y., 2018. Strong law of large numbers and Chover's law of the iterated logarithm under sub-linear expectations. Journal of Mathematical Analysis and Applications, 460(1), 252-270.

\bibitem{RefJ27} Wu Q.Y., Jiang Y.Y., 2020. Complete convergence and complete moment convergence for negatively dependent random variables under sub-linear expectations. Filomat, 34(4), 1093-1104.

\bibitem{RefJ28} Wu Y., Wang X.J., Shen A.T., 2021. Strong convergence properties for weighted sums of $m$-asymptotic negatively associated random variables and statistical applications. Statistical Papers, 62, 2169-2194.

\bibitem{RefJ29} Wu Y., Wang X.J., Zhang L.X., 2018. On the asymptotic approximation of inverse moment under sub-linear expectations. Journal of Mathematical Analysis and Applications, 468, 182-196.

\bibitem{RefJ30} Xi M.M., Wang X.J., Wu Y., 2019. Complete convergence for arrays of rowwise END random variables and its statistical applications under sub-linear expectations. Journal of the Korean Statistical Society, 48(3), 412-425.

\bibitem{RefJ31} Zhang L.X., 2015. Donsker's invariance principle under the sub-linear expectation with an application to Chung's law of the iterated logarithm. Communications in Mathematics and Statistics, 3(2), 187-214.

\bibitem{RefJ32} Zhang L.X., 2016a. Exponential inequalities under sub-linear expectations with applications to laws of the iterated logarithm. Science China Mathematics, 59(12), 2503-2526.

\bibitem{RefJ33} Zhang L.X., 2016b. Rosenthal's inequalities for independent and negatively dependent random variables under sub-linear expectations with applications. Science China Mathematics, 59(4), 751-768.

\bibitem{RefJ34} Zhang L.X., 2020. The convergence of the sums of independent random variables under the sub-linear expectations. Acta Mathematica Sinica, English Series, 36(3), 224-244.

\bibitem{RefJ35} Zhang L.X., 2022. Strong limit theorems for extended independent random variables and extended negatively dependent random variables under sub-linear expectations. Acta Mathematica Scientia, 42B(2), 467-490.

\bibitem{RefJ36} Zhong H.Y., Wu Q.Y., 2017. Complete convergence and complete moment convergence for weighted sums of extended negatively dependent random variables under sub-linear expectation. Journal of Inequalities and Applications, Volumn 2017, Article ID 261, 14 pages.

\end{thebibliography}
\end{document}